\documentclass[11pt]{article}
\usepackage{amsmath}
\usepackage{amsfonts}
\usepackage{amssymb}
\usepackage{hyperref}
\usepackage{indentfirst}
\usepackage{mathrsfs}

\usepackage{url}
\usepackage{color}

\usepackage[top=1in,bottom=1in,left=1.25in,right=1.25in]{geometry}

\date{}
\allowdisplaybreaks

\begin{document}

\renewcommand{\baselinestretch}{1.2}
\renewcommand{\arraystretch}{1.0}

\title{\bf Notes on multiplier Hopf algebras and invariants of framed links and 3-manifolds}
\author
{\textbf{Tao Yang} \footnote{College of Science, Nanjing Agricultural University, Nanjing 210095, China}
 $\cdot$
 \textbf{David Yetter} \footnote{Department of Mathematics, Kansas State University, Manhattan, KS 66506, US}
 }
\maketitle

\begin{center}
\begin{minipage}{12.cm}

 \textbf{Abstract}
 In this paper, we show that Hennings construction of invariants of framed links and 3-manifolds obtained from Hopf algebras 
 can also be carried out for some algebraic quantum groups.

\\


\end{minipage}
\end{center}
\normalsize

\section*{Introduction}
\def\theequation{\thesection.\arabic{equation}}
\setcounter{equation}{0}

 It is well known that Hopf algebras have been successfully used in the construction of invariants of links and 3-manifolds.  Several different constructions of invariants of links, framed links and 3-manifolds from data involving a Hopf algebra have been given.  Most involve using the Hopf algebra to produce a braided monoidal category by forming the category of modules, comodules or Yetter-Drinfel'd modules over the Hopf algebra.  One construction, due to Hennings \cite{H96} constructs invariants of framed links and 3-manifolds directly from a Hopf algebra equipped with a suitable trace function, without passing to any associated category.
 In his paper \cite{H96}, Hennings restricted his attention to the case of finite-dimensional Hopf algebras $H$, 
 because of the existence of a non-zero right integral in $H^*$. 
 In principle, his construction also holds for infinite-dimensional Hopf algebras $H$ whose dual possesses a non-zero right integral.  We refer to invariants arising from this construction as Hennings invariants.

 Multiplier Hopf algebras, introduced by A. Van Daele in \cite{V94}, are a non-unital generalization of Hopf algebras, 
 in which the target of the comultiplication is no longer the twofold tensor product of the underlying algebra, but an enlarged multiplier algebra. 
 A remarkable feature of multiplier Hopf algebras is that they admit a nice duality which extends the (incomplete) duality of Hopf algebras.
 This duality is based on left- and right-invariant linear functionals called integrals, which are analogues of the Haar measures of a group. 
 
 Hence, it is natural to wonder if multiplier Hopf algebras can be used in the invariant construction. 
 
 The primary purpose of this note is to generalize the technical results used in constructing Hennings invariants from Hopf algebras to the multiplier Hopf algebras. 
We find that the Hennings construction can be carried out for some instances of  the family of multiplier Hopf algebras know as algebraic quantum groups.  Along the way we also show that categories of finite dimensional modules over a slightly broader class of multiplier Hopf algebras -- ribbon quasitriangular multiplier Hopf algebras -- provide adequate data for construction of invariants of framed links.

 The paper is organized as follow:
 In Section 1, we recall some notions which we will use throughout, such as the definition of
 multiplier Hopf algebra, algebraic quantum groups and their duality.
 In section 2, we consider the properties of the ribbon element in a multiplier Hopf algebra, 
 and characterize ribbon quasitriangular multiplier Hopf algebras, also showing that
the category of finite-dimensional representations of a ribbon quasitriangular multiplier Hopf algebra is a ribbon category.
 In Section 3, we investigate the trace functions constructed from integrals on multiplier Hopf algebras and show that 
 Hennings construction can be carried out for suitable multiplier Hopf algebras.

\section{Preliminaries}
\def\theequation{\thesection.\arabic{equation}}
\setcounter{equation}{0}

 Throughout this paper, all vector spaces we considered are over a fixed field $k$ (such as the field $\mathbb{C}$ of the complex numbers).
 Let $A$ be an (associative) algebra. We do not assume that $A$ has a unit, but we do require that
 the product, seen as a bilinear form, is non-degenerated. This means that, whenever $a\in A$ and $ab=0$ for all $b\in A$
 or $ba=0$ for all $b\in A$, we must have that $a=0$.
 Then we can consider the multiplier algebra $M(A)$ of $A$.
 Recall that $M(A)$ is characterized as the largest algebra with identity containing $A$ as an essential two-sided ideal.  Elements of $M(A)$ are called multipliers on $A$ and we denote the identity multiplier by $1$.
 In particularly, we still have that, whenever $a\in M(A)$ and $ab=0$ for all $b\in A$ or $ba=0$ for all $b\in A$, again $a=0$.
 Furthermore, we consider the tensor algebra $A\otimes A$. It is again non-degenerated and we have its multiplier algebra $M(A\otimes A)$.
 There are natural imbeddings
 $$A\otimes A \subseteq M(A)\otimes M(A) \subseteq M(A\otimes A).$$
 In generally, when $A$ has no identity, these two inclusions are strict.
 If $A$ already has an identity, the product is obviously non-degenerate
 and $M(A)=A$ and $M(A\otimes A) = A\otimes A$. For more details about the concept of the multiplier algebra of an algebra, we refer the reader to \cite{V94}.

 Let $A$ and $B$ be non-degenerate algebras, if homomorphism $f: A\longrightarrow M(B)$ is non-degenerated
 (i.e., $f(A)B=B$ and $Bf(A)=B$),
 then has a unique extension to a homomorphism $M(A)\longrightarrow M(B)$, we also denote it $f$.

\subsection{Multiplier Hopf algebras and modules}

 Now, we recall the definition of a multiplier Hopf algebra (see \cite{V94} for details).
 A comultiplication on algebra $A$ is a homomorphism $\Delta: A \longrightarrow M(A \otimes A)$ such that $\Delta(a)(1 \otimes b)$ and
 $(a \otimes 1)\Delta(b)$ belong to $A\otimes A$ for all $a, b \in A$. We require $\Delta$ to be coassociative in the sense that
 \begin{eqnarray*}
 (a\otimes 1\otimes 1)(\Delta \otimes \iota)(\Delta(b)(1\otimes c))
 = (\iota \otimes \Delta)((a \otimes 1)\Delta(b))(1\otimes 1\otimes c)
 \end{eqnarray*}
 for all $a, b, c \in A$ (where $\iota$ denotes the identity map).

 A pair $(A, \Delta)$ of an algebra $A$ with non-degenerate product and a comultiplication $\Delta$ on $A$ is called
 a \emph{multiplier Hopf algebra}, if the linear map $T_{1}, T_{2}$ defined by
 \begin{eqnarray}
 T_{1}(a\otimes b)=\Delta(a)(1 \otimes b), \qquad T_{2}(a\otimes b)=(a \otimes 1)\Delta(b)
 \end{eqnarray}
 are bijective.

 The bijectivity of the above two maps is equivalent to the existence of a counit and an antipode S satisfying (and defined by)
 \begin{eqnarray}
 && (\varepsilon\otimes\iota)(\Delta(a)(1\otimes b)) = ab, \qquad  m(S\otimes\iota)(\Delta(a)(1\otimes b))=\varepsilon(a)b, \label{1.1} \\
 && (\iota\otimes\varepsilon)((a\otimes 1)\Delta(b)) = ab, \qquad  m(\iota\otimes S)((a\otimes 1)\Delta(b))=\varepsilon(b)a,\label{1.2}
 \end{eqnarray}
 where $\varepsilon:A\longrightarrow k$ is a homomorphism, $S: A\longrightarrow M(A)$ is an anti-homomorphism
 and $m$ is the multiplication map, considered as a linear map from $A\otimes A$ to $A$ and extended to
 $M(A)\otimes A$ and $A\otimes M(A)$.

 A multiplier Hopf algebra $(A, \Delta)$ is called \emph{regular} if $(A, \Delta^{cop})$ is also a multiplier Hopf algebra,
 where $\Delta^{cop}$ denotes the co-opposite comultiplication defined as $\Delta^{cop}=\tau \circ \Delta$ with $\tau$ the usual flip map
 from $A\otimes A$ to itself extended to $M(A\otimes A)$. In this case, $\Delta(a)(b \otimes 1), (1 \otimes a)\Delta(b) \in A \otimes A$
 for all $a, b\in A$.
 By Proposition 2.9 in \cite{V98}, multiplier Hopf algebra $(A, \Delta)$ is regular if and only if the antipode $S$ is bijective from $A$ to $A$.

 We will use the adapted Sweedler notation for multiplier Hopf algebras (see \cite{V08}), e.g., write $a_{(1)} \otimes a_{(2)}b$ for $\Delta(a)(1 \otimes b)$ 
 and $ab_{(1)} \otimes b_{(2)}$ for $(a \otimes 1)\Delta(b)$.
 \\
 
 Let $A$ be a regular multiplier Hopf algebra. Suppose $M$ is a left $A$-module with the module structure map $\cdot : A\otimes M\longrightarrow M$. 
 We will always assume that the module is non-degenerate, this means that $x=0$ if $x\in X$ and $a\cdot x=0$ for all $a\in A$.
 If the module is unital (i.e., $A\cdot M = M$), then it is non-degenerate and there exists an extension of the module structure to $M(A)$, this
 means that we can define $f\cdot m$, where $f\in M(A)$ and $m\in M$.
 In fact, since $m\in M = A\cdot M$, then $m= \sum_{i} a_{i}\cdot m_{i}$, $f\cdot m=\sum_{i} (fa_{i})\cdot m_{i}$.
 In this setting, we can easily get $1_{M(A)}\cdot x=x$. 
 For more details we refer the reader to \cite{DrVZ99}.

\subsection{Algebraic quantum groups and their dualities}

 Assume in what follows that $(A, \Delta)$ is a regular multiplier Hopf algebra. A linear functional $\varphi$ on $A$ is called left invariant if
 $(\iota \otimes \varphi)\Delta(a) = \varphi(a)1$ in $M(A)$ for all $a \in A$. A non-zero left invariant functional $\varphi$ is called a
 \emph{left integral} on $A$. A right integral $\psi$ can be defined similarly.

 In general, left and right integrals are unique up to a scalar if they exist.
 And if a left integral $\varphi$ exists, a right integral also exists, namely $\varphi \circ S$.  In general right and left integrals do not coincide.
 For a left integral $\varphi$, there is a unique group-like element called the {\em modular element} $\delta \in M(A)$ such that
 $\varphi(S(a)) = \varphi(a\delta)$ for all $a \in A$.

A regular multiplier Hopf algebra which admit a left integral is called an {\em algebraic quantum group} in \cite{V98}. 
 For an algebraic quantum group $(A, \Delta)$ with integral $\varphi$, define $\widehat{A}$ as the space of linear functionals on $A$
 of the form $\varphi(\cdot a)$, where $a \in A$. Then $\widehat{A}$ can be made into a regular multiplier Hopf algebra
 with a product(resp. coproduct $\widehat{\Delta}$ on $\widehat{A}$) dual to the coproduct $\Delta$ on $A$ (resp. product of $A$).
 It is called the dual of $(A, \Delta)$.
 The various objects associated with $(\widehat{A}, \widehat{\Delta})$ are denoted as for $(A, \Delta)$ but with a hat.
 However, we use $\varepsilon$ and $S$ also for the counit and antipode on the dual.
 The dual $(\widehat{A}, \widehat{\Delta})$ also has integrals, i.e., the dual is also an algebraic quantum group.
 A right integral $\widehat{\psi}$ on $\widehat{A}$ is defined by $\widehat{\psi}(\varphi(\cdot a)) = \varepsilon(a)$ for all $a\in A$.
 Repeating the procedure, i.e., taking the dual of $(\widehat{A}, \widehat{\Delta})$,
 we can get $\widehat{\widehat{A}}\cong A$ (see Theorem (Biduality) 4.12 in \cite{V98}).

 From Definition 1.5 in \cite{D10}, an algebraic quantum group $A$ is called \emph{counimodular},
 if the dual multiplier Hopf algebra $\widehat{A}$ is unimodular integral, i.e., $\widehat{\delta}= 1$ in $M(\widehat{A})$.
 For a counimodular algebraic quantum group, we have $\varphi(ab)=\varphi(bS^{2}(a))$ and $\psi(S^{2}(b)a) = \psi (a b)$ for
 all $a, b \in A$ (see Proposition 1.6 in \cite{D10}).

\section{Ribbon Quasitriangular Multiplier Hopf Algebras}
\def\theequation{\thesection.\arabic{equation}}
\setcounter{equation}{0}

 Recall from \cite{DVW05}, that a regular multiplier Hopf algebra is quasitriangular if it is equipped with an invertible multiplier $R \in M(A\otimes A)$ satisfying

\begin{eqnarray} 
R\Delta(a) & = & \Delta^{cop}(a)R \\
(\Delta \otimes \iota)(R) & = & R^{13}R^{23} \\
(\iota \otimes \Delta)(R) & = & R^{13}R^{12} 
\end{eqnarray}

\noindent Such a multiplier $R$ is called a generalized $R$-matrix for $A$.
 
As in \cite{DVW05}, if $R$ is a  generalized $R$-matrix for a regular multiplier Hopf algebra $A$,
 i.e., furthermore $R(a\otimes 1), (a\otimes 1)R \in A\otimes M(A)$ for all $a\in A$.
 Then we can define a multiplier $u$ such that for all $a\in A$, 
 \begin{eqnarray}
 ua=S(R^{(2)})R^{(1)}a.
 \end{eqnarray}
 Because of the invertibility of $R$, by Lemma 2.2 and corollary 2.5 in \cite{DVW05}, we have
 $u$ is invertible in $M(A)$ and satisfies the following identities: for any $a\in A$,
 \begin{eqnarray}
 && u^{-1}a = S^{-2}(R^{(2)})R^{(1)}a = S^{-1}(R^{(2)})S(R^{(1)})a = R^{(2)} S^{2}(R^{(1)})a, \\
 && S^{2}(a) = uau^{-1} = S(u)^{-1}aS(u), \\
 && \Delta(u) = (R_{21}R_{12})^{-1}(u\otimes u) = (u\otimes u)(R_{21}R_{12})^{-1} \in M(A\otimes A).
 \end{eqnarray}

 \textbf{Proposition \thesection.1}
 The multiplier $uS(u)$ belongs to $Z(M(A))$, the center of $M(A)$.

 \emph{Proof} Indeed, for any $a\in A$,
 \begin{eqnarray*}
 S^{2}(a) &=& u a u^{-1} = u S(u) \underline{S(u)^{-1}aS(u)} S(u)^{-1} u^{-1} \\
 &=& u S(u) S^{2}(a) S(u)^{-1} u^{-1}.
 \end{eqnarray*}
 Thus, because $S^{2}(A)=A$ we have that $a u S(u) = u S(u) a$ for any $a \in A$.

 For any $x\in M(A)$,
 \begin{eqnarray*}
 a(x u S(u)) &=& (ax) u S(u) = u S(u)(ax) \\
 &=& (u S(u)a) x = (a u S(u)) x \\
 &=& a (u S(u) x)
 \end{eqnarray*}
 holds for all $a\in A$, thus $x u S(u) = u S(u) x$.
 $\hfill \Box$
 \\

 \textbf{Remark} The second part of the proof actually follows the fact that
 \begin{eqnarray*}
 \{z\in M(A)| za=az \mbox{ for all } a\in A\}
 = \{z\in M(A)| zx=xz \mbox{ for all } x\in M(A)\}.
 \end{eqnarray*}
Plainly $Z(A)$, the center of $A$ is contained in $Z(M(A))$, the center of the multiplier algebra of $A$.
 \\

 \textbf{Definition \thesection.2} (\cite{I15})
 Let $A$ be a quasitriangular multiplier Hopf algebra with generalized $R$-matrix $R$, $(A, R)$ is called \emph{ribbon multiplier Hopf algebra}
 if there exists a central element $v$ in $Z(M(A))$ such that:
 \begin{eqnarray*}
 && v^{2} = u S(u), \quad S(v)=v, \quad \varepsilon(v)=1, \\
 && \Delta(v) = (R_{21}R_{12})^{-1}(v\otimes v).
 \end{eqnarray*}

 \textbf{Remark}
 (1) This ribbon structure, introduced in Definition 6.4 of \cite{I15}, 
 is a special case of Definition 2.2 in \cite{D10}, in which $G$ is the trivial group.

 (2) Since $u$ is invertible, $v$ (called the \emph{ribbon element}) is also invertible in $M(A)$.
 The multiplier $u v^{-1}$ is group-like, and so $S(u v^{-1}) = v u^{-1}$.

 (3) If $v$ exists it may be non-unique. In order to describe the non-uniqueness of $v$,
 we define the abelian group $E(M(A))\subset M(A)$ consisting of central elements $E\in M(A)$ such that
 \begin{eqnarray*}
 E^{2}=1, \quad S(E)=E, \quad \varepsilon(E)=1, \quad \Delta(E)=E\otimes E.
 \end{eqnarray*}
 It is obvious that for any ribbon quasitriangular multiplier Hopf algebra $(A, R, v)$ and any $E\in M(A)$ the triple
 $(A, R, Ev)$ is ribbon  iff $E \in E(M(A))$.
 \\

 \textbf{Example \thesection.3}
 (1) Let $A$ be an infinite-dimensional co-Frobenius Hopf algebra, 
 and suppose that $A$ is coquasitriangular (or cobraided) in the sense of Definition VIII.5.1 in \cite{K95}.
 Then the dual $\widehat{A}$ is a discrete quasitriangular multiplier Hopf algebra rather than a classical Hopf algebra.
 Furthermore, if $A$ is coribbon, then by Lemma A.2 in \cite{D10}, $\widehat{A}$ is a ribbon quasitriangular multiplier Hopf algebra.
 
 (2) By Example A.3 in \cite{D10} and Proposition 6.16 in \cite{I15}, $D_{q}(sl_{2})$ and 
 $\mathcal{\hat{U}}_{q\tilde{q}}(\mathfrak{g}_{\mathbb{R}})$ are ribbon multiplier Hopf algebras.
 
  (3) For an infinite group $G$, denote by $A=kG$ the corresponding group algebra 
 and by $\widehat{A} = k(G)$ the classical motivation example introduced in \cite{V94}. 
 Let $D(G) = A \otimes \widehat{A}$ as a vector space and equip it with a new multiplication
 \begin{eqnarray*}
 (g_{1} \delta_{h_{1}}) \cdot (g_{2} \delta_{h_{2}}) = g_{1} g_{2} \delta_{g_{2}^{-1} h_{1} g_{2} }\delta_{h_{2}},
 \end{eqnarray*}
 where we have dropped the $\otimes$-symbol on both sides, and the coproduct, counit, antipode  are given as follows:
 for any $g \delta_{h} \in D(G)$,
 \begin{eqnarray*}
 \Delta(g \delta_{h}) &=& \Delta_{A}(g)\Delta_{\widehat{A}}^{cop} (\delta_{h}) = \sum_{p\in G} g\delta_{p} \otimes g\delta_{hp^{-1}}, \\
 \varepsilon (g \delta_{h}) &=& \varepsilon_{A} (g) \varepsilon_{\widehat{A}} (\delta_{h}) = \delta_{h, e}, \\
 S(g \delta_{h}) &=& S_{\widehat{A}}(\delta_{h}) S_{A}(g) = \delta_{h^{-1}} g^{-1} = g^{-1} \delta_{g h^{-1} g^{-1}}.
 \end{eqnarray*}
 Then follow the computation in \cite{P11}
 Drinfel'd double $D(G)$ is a ribbon quasitriangular multiplier Hopf algebra with generalized $R$-matrix $R$ and ribbon element $v$ as follows
 \begin{eqnarray*}
 R &=& \sum_{g\in G} g\otimes \delta_{g} \in M(D(G) \otimes D(G)), \\
 v &=& \sum_{g\in G} g^{-1} \delta_{g} \in M(D(G)).
 \end{eqnarray*}

The following theorem characterizes which quasitriangular multiplier Hopf algebras are ribbon.
The proof uses the next lemma.

 \textbf{Lemma \thesection.4}
 (1) For any multiplier $x\in M(A)$, $S^{2}(x) = uxu^{-1} = S(u)^{-1}xS(u)$.

 (2) For any group-like elements $\mathrm{g}\in M(A)$, $u\mathrm{g} = \mathrm{g}u$.

 \emph{Proof} (1) We only check $S^{2}(x) = uxu^{-1}$, the other part is similar.
 For any $a\in A$, $xa\in A$ and $S^{2}(x)S^{2}(a) = S^{2}(xa) = uxau^{-1} = uxu^{-1}uau^{-1} = uxu^{-1}S^{2}(a)$.
 Since $S$ is bijective, $S^{2}(x)a = uxu^{-1}a$ holds for any $a$ in $A$, hence the result holds.

 (2) This follows easily from $\mathrm{g} = S^{2}(\mathrm{g}) = u\mathrm{g}u^{-1}$.
 $\hfill \Box$
 \\

 \textbf{Theorem \thesection.5}
 If $(A, R)$ is a quasitriangular multiplier Hopf algebra, then $A$ is a ribbon multiplier Hopf algebra
 if and only if there exists a group-like element $g\in M(A)$ such that $g^{2} = S(u)^{-1} u$ and $S^{2}(a) = gag^{-1}$ for all $a\in A$.

 \emph{Proof}
 Let $A$ be a ribbon multiplier Hopf algebra with ribbon element $v$, then $S^{2}(a) = uau^{-1} = uv^{-1}vau^{-1} = uv^{-1}avu^{-1}$.
 Set $g=uv^{-1}$, then $S^{2}(a) = gag^{-1}$,
 \begin{eqnarray*}
 && \Delta(g) = \Delta(uv^{-1}) = \Delta(u)\Delta(v^{-1})  = (u\otimes u)(R_{21}R_{12})^{-1} (v^{-1}\otimes v^{-1})(R_{21}R_{12}) = g\otimes g, \\
 && g^{2} = (uv^{-1})^{2} = v^{-2}u^{2} = (uS(u))^{-1}u^{2} = S(u)^{-1}u.
 \end{eqnarray*}

 Conversely,  if there exists an element $g\in M(A)$ such that $g^{2} = S(u)^{-1} u$ and $S^{2}(a) = gag^{-1}$,
 let $v=ug^{-1}$, then we can easily check $v$ satisfies the conditions of Definition \thesection.2:
 \begin{eqnarray*}
 && v^{2} =  (ug^{-1})^{2} = u^{2}g^{-2} = u^{2}u^{-1}S(u) = uS(u), \\
 && S(v) = S(ug^{-1}) = gS(u) = g (u g^{-2}) = ug^{-1} = v, \\
 && \Delta(v) = \Delta(ug^{-1}) = \Delta(u)\Delta(g^{-1}) = (R_{21}R_{12})^{-1}(u\otimes u) (g^{-1}\otimes g^{-1}) = (R_{21}R_{12})^{-1}(v\otimes v).
 \end{eqnarray*}
To see that $v$ is central observe for any $a\in A$,
 $va=av \Leftrightarrow g^{-1}ua = ag^{-1}u \Leftrightarrow uau^{-1} = gag^{-1}$, but both of sides
of the last equation equal $S^2(a)$, by Lemma \thesection.4 and by hypothesis, respectively.
 $\hfill \Box$
 \\

 \textbf{Remark} (1) $S(u)^{-1} u = u S(u)^{-1}$. In fact, from $g^{2} = S(u)^{-1} u$ we get $S(u)^{-1} = g^{2}u^{-1}$,
 so $u S(u)^{-1} = ug^{2}u^{-1} =g^{2}$. The last equation holds by Lemma 3.4 (2).
 
 (2) This result generalizes Proposition 4.1 in \cite{H96} and Theorem 2 (2) in \cite{L10}. 
 This theorem also implies that there is a one-to-one correspondence between the ribbon elements and the group-like
 elements satisfying certain conditions. 
 \\
 
 \textbf{Corollary \thesection.6}
 Let $(A, R)$ be a quasitriangular multiplier Hopf algebra in which the order of $G(M(A))$, the group of group-like multipliers, is odd. Then
 $(A, R)$ is ribbon if and only if $S^{2|G(M(A))|} = \iota$.
 
 \emph{Proof}
 If $A$ is ribbon with ribbon element $v$, let $g=uv^{-1}$ as above, then $S^{2}(a) = gag^{-1}$
 and $S^{2|G(M(A))|}(a) = g^{|G(M(A))|} a g^{-|G(M(A))|}$ = a, i.e., $S^{2|G(M(A))|} = \iota$.
 
 Conversely, if $S^{2|G(M(A))|} = \iota$, writing $|G(M(A))| = 2n+1$ we have that $(u S(u)^{-1})^{2n+2} = u S(u)^{-1}$ because of $u S(u)^{-1}\in G(M(A))$. 
 Let $g = (u S(u)^{-1})^{n+1}$, then $g$ is a group-like element in $M(A)$ and $g^{2} = u S(u)^{-1} = S(u)^{-1} u$, while
 \begin{eqnarray*}
 gag^{-1} 
 &=& (u S(u)^{-1})^{n+1} a (u S(u)^{-1})^{-(n+1)} 
  = u^{n+1} S(u)^{-(n+1)} a S(u)^{n+1} u^{-(n+1)} \\
 &=& u^{n+1} S^{2(n+1)}(a) u^{-(n+1)}
 = S^{4(n+1)}(a) = S^{2|G(M(A))|+2}(a) \\
 &=& S^{2}(a).
 \end{eqnarray*}
 Following Theorem \thesection.5, we can get $A$ is a ribbon multiplier Hopf algebra.
 $\hfill \Box$
 \\
 
 Let  ${}_{A}\mathcal{M}$ be the category of representations of a multiplier Hopf algebra $A$, 
 whose objects are unital left $A$-modules and whose morphisms are $A$-linear homomorphisms.
 For any unital left $A$-module $M$ and $N$, we can easily check that $M\otimes N$ is also a unital left $A$-module with the action
 \begin{eqnarray}
 a\cdot (m\otimes n) = a_{(1)}\cdot m \otimes a_{(2)}\cdot n,
 \end{eqnarray}
 where $a\in A$, $m\in M$ and $n\in N$. 
 This module action of $A$ on $M\otimes N$ makes sense, since $M$ and $N$ are unital $A$-modules, 
 it is clear that $M\otimes N$ is an $A\otimes A$-module by $(a\otimes a')\cdot (m\otimes n) = a\cdot m \otimes a'\cdot n$
 and by Proposition 3.3 in \cite{DrVZ99} we can extend it to $M(A\otimes A)$. 
 Also the action of $A$ on $M\otimes N$ is unital.


\noindent   By Proposition 3.4 in \cite{DrVZ99}, ${}_{A}\mathcal{M}$ is a monoidal category.
 The unit object is the basic field $k$ with the module action $a\cdot \kappa = \varepsilon(a)\kappa$ and the unit constraints are given by
  \begin{eqnarray*}
 l_{M}: k\otimes M\rightarrow M, && \kappa\otimes m\mapsto \kappa m, \\
 r_{M}: M\otimes k\rightarrow M, && m\otimes\kappa\mapsto \kappa m.
 \end{eqnarray*} 
 If, in addition, $A$ is a quasitriangular multiplier Hopf algebra, then by Theorem 4 in \cite{Z99} ${}_{A}\mathcal{M}$ is a braided monoidal category
 with the braiding given by 
 \begin{eqnarray*}
 c_{M, N} (m\otimes n) = R^{(2)}\cdot n \otimes R^{(1)}\cdot m.
 \end{eqnarray*}
 
 In the following, we consider the full subcategory ${}_{A}\mathcal{M}_{f}$, in which all objects are finite-dimensional and morphisms are $A$-linear homomorphisms. 
 Similar to Proposition 3.1 in \cite{YZM13}, we can get that ${}_{A}\mathcal{M}_{f}$ is a tensor category with left and right duality as follows:
 For any finite-dimensional unital left $A$-module $M$, $M^* = Hom(M, k)$ becomes an object in ${}_{A}\mathcal{M}_{f}$ with module action
 \begin{eqnarray}
 (a\cdot f)(m) = f(S(a)\cdot m), \quad \forall a\in A, m\in M, f\in M^*. \label{Y1}
 \end{eqnarray}
 Moreover, the maps $b_{M}: k\rightarrow M\otimes M^*$, $b_{M}(1_{k})=\sum_{i} e_{i} \otimes e^{i}$ and
 $d_{M}: M^*\otimes M\rightarrow k$, $d_{M}(f\otimes m)=f(m)$ are morphisms in ${}_{A}\mathcal{M}_{f}$, where
 $\{ e_{i} \}$ and $\{ e^{i} \}$ are dual bases in $M$ and $M^*$, and $(id_{M}\otimes d_{M})(b_{M}\otimes id_{M}) = id_{M}$,
 $(d_{M} \otimes id_{M^*})(id_{M^*}\otimes b_{M})=id_{M^*}$.

 Similarly, ${}^*M = Hom(M, k)$ becomes an object in ${}_{A}\mathcal{M}_{f}$ with module action
 \begin{eqnarray*}
 (a\cdot f)(m) = f(S^{-1}(a) \cdot m), \quad \forall a\in A, m\in M, f\in {}^*M.
 \end{eqnarray*}
 Moreover, the maps $b_{M}: k\rightarrow {}^*M\otimes M$, $b_{M}(1_{k})=\sum_{i} e^{i} \otimes e_{i}$ and
 $d_{M}: M\otimes {}^*M\rightarrow k$, $d_{M}(m\otimes f)=f(m)$ are morphisms in ${}_{A}\mathcal{M}_{f}$,
 and $(d_{M}\otimes id_{M})(id_{M}\otimes b_{M}) = id_{M}$, $(id_{{}^*M} \otimes d_{M})(b_{M}\otimes id_{{}^*M})=id_{{}^*M}$.
 Following these results, we can get that
 \\
 
 \textbf{Proposition \thesection.7} The category ${}_{A}\mathcal{M}_{f}$ is rigid.  
 \\
 
 If, furthermore, $A$ is a ribbon quasitriangular multiplier Hopf algebra, then each $M\in {}_{A}\mathcal{M}_{f}$ has a twist map
 \begin{eqnarray*}
 \theta_{M}: M \longrightarrow M, \quad m \mapsto v^{-1}\cdot m.
 \end{eqnarray*}
 Then we can get the following result:
 \\
 
 \textbf{Theorem \thesection.8}
 For a ribbon quasitriangular multiplier Hopf algebra $A$, the category ${}_{A}\mathcal{M}_{f}$ is a ribbon category.
 
 \emph{Proof}
 For any morphism $f: M\rightarrow N$ in ${}_{A}\mathcal{M}_{f}$, we can easily get that $f\circ \theta_{M} = \theta_{N} \circ f$ because $f$ is an $A$-module morphism.
 In the following, we need to check that $\theta_{M\otimes N} = (\theta_{M}\otimes \theta_{N})c_{N, M}c_{M, N}$ and $(\theta_{M})^*=\theta_{M^*}$.
 Indeed, for all $a\in A$, $m\in M$ and $n\in N$,
 \begin{eqnarray*}
 (\theta_{M}\otimes \theta_{N})c_{N, M}c_{M, N} (m\otimes n)
 &=& (\theta_{M}\otimes \theta_{N})\Big( R_{21} R_{12} (m\otimes n) \Big) \\
 &=& (v^{-1}\otimes v^{-1}) R_{21} R_{12} (m\otimes n) \\
 &=& \Delta(v^{-1}) (m\otimes n) \\
 &=& \theta_{M\otimes N} (m\otimes n).
 \end{eqnarray*}
 And for any $a\in A$, $m\in M$ and $f\in M^*$,
  \begin{eqnarray*}
 (\theta_{M})^*(f)(m)
 &=& f\Big(\theta_{M}(m)\Big)
 = f\Big(v^{-1}\cdot m\Big) \\
 &=& f\Big(S(v^{-1})\cdot m\Big)
 = (v^{-1}\cdot f)(m) \\
 &=& \theta_{M^*}(f)(m).
 \end{eqnarray*} 
 \noindent where the next to last equality holds by \ref{Y1}, above.  Following Definition XIV.3.2 in \cite{K95},  ${}_{A}\mathcal{M}_{f}$ is a ribbon category.
 $\hfill \Box$
 \\
 
It follows from the coherence theorem of Shum \cite{S} that a choice of object (resp. of $n$ objects) in ${}_{A}\mathcal{M}_{f}$ gives rise to an invariant of framed links (resp. $n$-colored framed links).   All these invariants can also be obtained from categories of comodules over ordinary Hopf algebras, since the ribbon subcategory of ${}_{A}\mathcal{M}_{f}$ generated by the chosen object(s) admits a fiber functor to the category of vector spaces and will, thus be equivalent to a subcategory of the category of comodules of the ribbon Hopf algebra obtained by Tannaka reconstruction.  It would seem to be the case that the construction of the invariant from a multiplier Hopf algebra by way of its category of finite dimensional representations is more natural than construction from an ordinary Hopf algebra in cases where the latter can be found only retrospectively via Tannaka reconstruction after the invariant is already in hand.

\section{Trace functions and the Hennings construction}
\def\theequation{\thesection.\arabic{equation}}
\setcounter{equation}{0}

In the Hennings \cite{H96}, the ribbon element and trace functions play an essential role in the construction of invariants of links and 3-manifolds.
Above, we have considered the ribbon element. In this section, we will consider the trace functions on multiplier Hopf algebra 
 and generalize the Hennings construction from finite-dimensional Hopf algebras to algebraic quantum groups.
 \\

 Let $A$ be a regular multiplier Hopf algebra with right integral $\psi$, namely an algebraic quantum group. 
 Assume that $A$ is counimodular in the sense of Definition 1.5 in \cite{D10}, i.e., $\widehat{\delta} = 1$ in $M(A)$.
 Then by Proposition 1.6 in \cite{D10} we have
 \begin{eqnarray}
 \psi(S^{2}(b)a) = \psi (a b) \quad \mbox{for all} \quad a, b \in A. \label{3.1}
 \end{eqnarray}

 Now let
 \begin{eqnarray*}
 && C_{0} = \{\mu=\psi(m\cdot), m\in M(A) :  \mu(ab)=\mu(ba), \forall a, b\in A\}. 
 \end{eqnarray*}
 In the following, we concentrate on finding trace functions.  
 The following result shows a necessary and sufficient condition of the form for some trace functions when $A$ is counimodular,
 generalizing Proposition 4.2 of \cite{H96} to multiplier Hopf algebras:
 \\

 \textbf{Proposition \thesection.1}
 Let $(A, R)$ be a counimodular quasitriangular ribbon multiplier Hopf algebra with ribbon element $v$ and set $g=uv^{-1}$. 
 Then we have
 \begin{enumerate}
 \item[(1)] $\mu\in C_{0}$ if and only if $\mu=\mu_{z}=\psi(gz\cdot)$ for some $z\in Z(M(A))$,
 \item[(2)] $\mu\in C_{0}$ and $\mu = \mu \circ S$ if and only if $\mu=\mu_{z}$ for some $z = S(z)\in Z(M(A))$.
 \end{enumerate}

 \emph{Proof}
 (1) If $z\in Z(M(A))$, then for any $a, b\in A$, 
 \begin{eqnarray*}
 \mu_{z}(ab)=\psi(gzab)\stackrel{(\ref{3.1})}{=} \psi(S^{2}(b)gza)=\psi(gbza)=\psi(gzba)=\mu_{z}(ba),
 \end{eqnarray*}
 so $\mu\in C_{0}$.

 Conversely, if $z\in M(A)$ such that $\mu_{z}\in C_{0}$, then for any $a, b\in A$,
 \begin{eqnarray*}
 \psi(gzab)=\mu_{z}(ab)=\mu_{z}(ba)=\psi(gzba) \stackrel{(\ref{3.1})}{=} \psi(S^{2}(a)gzb)=\psi(gazb), 
 \end{eqnarray*}
 hence $\psi((gza-gaz)b)=0$ for any $b\in A$.
 From Proposition 3.4 in \cite{V98} non-zero invariant functionals are always faithful, so $gza=gaz$ for all $a\in A$,
 we can get $z\in Z(M(A))$.

 (2) By Proposition 1.6 (4) in \cite{D10}, $\psi(gzS(\cdot))=\psi(gz\cdot)$ if and only if $S(gz)=gz\delta$,
 i.e., $S(z) = gz\delta g = zg\delta g$, where $\delta$ is the modular element.
 From Proposition 2.9 in \cite{DVW05}, $uS(u)^{-1} = \delta^{-1}((\iota\otimes \langle\widehat{\delta}^{-1},\cdot\rangle)R)$.
 When $A$ is counimodular $\widehat{\delta}=1$, so $g^{2} = uS(u)^{-1} = \delta^{-1}$.
 Hence, $\psi(gzS(\cdot))=\psi(gz\cdot)$ if and only if $S(z)=z$.
 $\hfill \Box$
 \\


 Thus, in order to find an invariant of regular isotopy for unoriented framed links, as in the Hennings's construction,
 we need a counimodular ribbon quasitriangular multiplier Hopf algebra, together with a choice of a central mulitplier $z \in Z(M(A))$ satisfying $z = S(z)$.
 Such multiplier Hopf algebras exist, as, for example, $D_{q}(sl_{2})$ introduced in \cite{D10}.

In order to extend Hennings's construction of 3-manifold invariants, however, we will see that $z$ must lie in the center of the algebra, $Z(A)$, rather than having the freedom to lie in all of $Z(M(A))$.

 To proceed, we first recall from \cite{D10} a right $\widehat{A}$-module structure on $A$, denoted by $A\blacktriangleleft \widehat{A}$,
 which is a special case of $A\blacktriangleleft B$ under a multiplier Hopf algebra paring $\langle A, B \rangle$ introduced in \cite{DV04}.
 For $a, b\in A$, we define $a \blacktriangleleft \psi(b\cdot) = \psi(b a_{(1)}) a_{(2)} = (\psi\otimes\iota)((b\otimes 1)\Delta(a))$ in $A$.
 Observe that $A\blacktriangleleft \widehat{A}$ defines an $\widehat{A}$-module structure on $A$ because the product in $\widehat{A}$ 
 is dual to the coproduct in $A$. As $(A\otimes 1)\Delta(A)=A\otimes A$, we have $A\blacktriangleleft \widehat{A} = A$.

 In the following, as usual, $(A, R)$ is a counimodular quasitriangular ribbon multiplier Hopf algebra with ribbon element $v\in M(A)$.
 Suppose there is an element $z = S(z) \in Z(A)$, then $\mu_{z} \in \widehat{A}$ and $\mu_{z} = \mu_{z} \circ S = S\mu_{z}$.
 Consider the multipliers
 \begin{eqnarray*}
 X_{z} & =  & S(u^{-1})(S(u)\blacktriangleleft \mu_{z}) - \psi(zv)v^{-1} \;\mbox{\rm and} \\ Y_{z} & = & u(u^{-1}\blacktriangleleft \mu_{z}) - \psi(zv^{-1})v,
 \end{eqnarray*}
 where the multiplier $S(u)\blacktriangleleft \mu_{z} = (\psi\otimes \iota)((gz\otimes 1)(S(u)\otimes S(u))(R_{21}R_{12})^{-1})$ and
 $u^{-1}\blacktriangleleft \mu_{z} = (\psi\otimes \iota)((gz\otimes 1)(u^{-1}\otimes u^{-1}) R_{21}R_{12})$ belong to $M(A)$.
 \\

 \textbf{Remark} 
 (1) Here we require $z$ in $Z(A)$ rather than $Z(M(A))$ mainly because $S(u)\blacktriangleleft \mu_{z}$ is in general not well-defined when $z\in Z(M(A))$.
 This condition lets $X_{z}$ and $Y_{z}$ make sense in the framework of multiplier Hopf algebras.
 
 (2) If furthermore $A$ is of discrete type, i.e. $A$ has a left cointegral $t$. We can normalize it to satisfy $\psi(t)=1$.
 Since $A$ is counimodular, by Lemma 1.7 (1) in \cite{D10} $t$ is also a right integral, 
 and one example occurs by putting $z=t$, which implies the additional condition $z\in Z(A)$ is reasonable. 
 
 (3) An (infinite-dimensional) coFrobenius Hopf algebra $H$ can be treated as a trivial example of multiplier Hopf algebra. 
 In this case, $M(H)=H$ and $Z(M(H)) = Z(H)$ and $z=1$ trivially satisfied the required conditions. 
 \\

The following generalizes Proposition 4.3 of \cite{H96} to multiplier Hopf algebras, with (\ref{a0}) being Hennings' condition (1) and (\ref{a1}) and (\ref{a2}) being Hennings' condition (2):

 \textbf{Proposition \thesection.2}
 If element $z = S(z)\in Z(A)$ satisfying $(1\otimes z)\Delta(z) = z\otimes z$, then 
 \begin{eqnarray}
 zX_{z} = zY_{z} = 0,  \label{a0}
 \end{eqnarray}
 and in addition,
 \begin{eqnarray}
 && (z\otimes 1 \otimes \ldots \otimes 1)\Delta^{(n)}(vX_{z}) = z\otimes \Delta^{(n-1)}(vX_{z}), \label{a1}\\
 && (z\otimes 1 \otimes \ldots \otimes 1)\Delta^{(n)}(v^{-1}Y_{z}) = z\otimes \Delta^{(n-1)}(v^{-1}X_{z}), \label{a2}
 \end{eqnarray}
 where $\Delta^{(n)}: A \longrightarrow M(\otimes^{n} A)$ is linear map obtained from $n-1$ applications of the coproduct.

 \emph{Proof}
 First we compute $zX_{z}$ and $zY_{z}$ as follows.
 \begin{eqnarray*}
 zX_{z}
 &=& z S(u^{-1})(S(u)\blacktriangleleft \mu_{z}) - \psi(zv)zv^{-1} \\
 &=& z S(u^{-1})(\psi\otimes \iota)((gz\otimes 1)(S(u)\otimes S(u))(R_{21}R_{12})^{-1}) - \psi(zv)zv^{-1} \\
 &=& (\psi\otimes \iota)((gz S(u)\otimes z)(R_{21}R_{12})^{-1}) - \psi(zv)zv^{-1} \\
 &=& (\psi\otimes \iota)((vz\otimes z)(R_{21}R_{12})^{-1}) - \psi(zv)zv^{-1} \\
 &=& (\psi\otimes \iota)((1\otimes zv^{-1})\Delta(zv)) - \psi(zv)zv^{-1} \\
 &=& 0,
 \end{eqnarray*}
 and
 \begin{eqnarray*}
 zY_{z}
 &=& z u(u^{-1}\blacktriangleleft \mu_{z}) - \psi(zv^{-1})v \\
 &=& zu (\mu_{z} \otimes \iota)\Delta(u^{-1}) - \psi(zv^{-1})v \\
 &=& (\mu_{z} \otimes \iota)((1\otimes zu)\Delta(u^{-1})) - \psi(zv^{-1})v \\
 &=& (\psi \otimes \iota)((gz\otimes zu)\Delta(u^{-1})) - \psi(zv^{-1})v \\
 &=& (\psi \otimes \iota)((gz u^{-1}\otimes z)R_{21}R_{12}) - \psi(zv^{-1})v \\
 &=& v(\psi \otimes \iota)((z v^{-1}\otimes zv^{-1})R_{21}R_{12}) - \psi(zv^{-1})v \\
 &=& v(\psi \otimes \iota)\Delta(zv^{-1}) - \psi(zv^{-1})v \\
 &=& 0.
 \end{eqnarray*}

 Then we check the equation (\ref{a1}), and (\ref{a2}) is similar. Now we see
 $vX_{z} = v S(u^{-1})(S(u)\blacktriangleleft \mu_{z}) - \psi(zv)1_{M(A)}
 = g(S(u)\blacktriangleleft \mu_{z}) - \psi(zv)1_{M(A)}
 = \psi(gzS(u)_{(1)}) gS(u)_{(2)} - \psi(zv)1_{M(A)}$,
 and $\Delta^{(n)}: A \longrightarrow M(\otimes^{n} A)$ is non-degenerated, has a unique extension from $M(A)$ to $M(\otimes^{n} A)$.
 Then for any $a_{1}, a_{2}, ..., a_{n} \in A$,
 $\Delta^{(n)}(\psi(gzS(u)_{(1)}) gS(u)_{(2)})(a_{1}\otimes a_{2} \otimes ...\otimes a_{n}) 
 = \psi(gzS(u)_{(1)}) gS(u)_{(2)}a_{1} \otimes \ldots \otimes gS(u)_{(n+1)}a_{n} \in \otimes^{n} A$,
 so
 \begin{eqnarray*}
 && (z\otimes 1 \otimes \ldots \otimes 1)\Delta^{(n)}(vX_{z}) (a_{1}\otimes a_{2} \otimes ...\otimes a_{n}) \\
 &=& (z\otimes 1 \otimes \ldots \otimes 1)\Delta^{(n)}(\psi(gzS(u)_{(1)}) gS(u)_{(2)} - \psi(zv)1_{M(A)}) (a_{1}\otimes a_{2} \otimes ...\otimes a_{n}) \\
 &=& \psi(gzS(u)_{(1)}) gzS(u)_{(2)} a_{1} \otimes \ldots \otimes gS(u)_{(n+1)} a_{n}
    - \psi(zv)z a_{1}\otimes a_{2} \otimes\ldots \otimes a_{n} \\
 &=& (\psi\otimes \iota) \big((1\otimes z)\Delta(zgS(u)_{(1)}) \big) a_{1} \otimes \ldots \otimes gS(u)_{(n)} a_{n}
    - \psi(zv)z a_{1}\otimes a_{2} \otimes\ldots \otimes a_{n} \\
 &=& \psi(zgS(u)_{(1)}))z a_{1} \otimes \ldots \otimes gS(u)_{(n)} a_{n}
    - \psi(zv)z a_{1}\otimes a_{2} \otimes\ldots \otimes a_{n} \\
 &=& z a_{1}\otimes \Big( \psi(zgS(u)_{(1)})gS(u)_{(2)} a_{2} \otimes \ldots \otimes gS(u)_{(n)} a_{n}
    - \psi(zv) a_{2} \otimes\ldots \otimes a_{n} \Big) \\
 &=& \Big( z\otimes \big( \psi(zgS(u)_{(1)})gS(u)_{(2)} \otimes \ldots \otimes gS(u)_{(n)}
    - \psi(zv)1_{M(A)} \otimes\ldots \otimes 1_{M(A)} \big) \Big) \\
  && (a_{1}\otimes a_{2} \otimes ...\otimes a_{n}) \\
 &=& (z\otimes \Delta^{(n-1)}(vX_{z})) (a_{1}\otimes a_{2} \otimes ...\otimes a_{n}).
 \end{eqnarray*}
 Therefore equation (\ref{a1}) holds. This completes the proof.
 $\hfill \Box$
 \\

 It turns out that as in Hopf algebra case conditions (\ref{a0}), (\ref{a1}) and (\ref{a2}) are precisely what is needed to obtain framed link invariants
 which behave well under the Fenn-Rourke moves, and so which can be used to define invariants of 3-manifolds.
 We find that Hennings construction can also be carried out with multiplier Hopf algebras.
 It remains to provide an example showing that something is gained from the greater generality.  
 \\

 \textbf{Example \thesection.3} 
 Fix an infinite group $G$ with at least one non-trivial finite subgroup $K$, and consider the quantum double multiplier Hopf algebra of Example 2.3, $D(G)$.  
 First, $D(G)$ is unimodular with integral $\psi$ given on the basis $\{ g\delta_h | g,h \in G\}$ by $\psi(g\delta_h) = \delta_{g,e}$ 
 (where the second $\delta$ is the Kronecker $\delta$).  Thus, we wish to find elements $z \in D(G)$ satisfying  $z = S(z)$ and $(1\otimes z) \Delta(z) = z\otimes z$.

A brief consideration of the formulas for multiplication, comultiplication and antipode suggest that it would be reasonable to look for elements which are sums of basis elements $g\delta_h$, rather than more general linear combinations.  Let $\Sigma \subset G\times G$ be a finite subset, and consider the element 

\[ z = z_\Sigma := \sum_{(g,h)\in \Sigma} g\delta_h \]

Then $S(z) = \sum_{(g,h) \in \Sigma} g^{-1}\delta_{gh^{-1}g^{-1}}$, so the condition $z = S(z)$ will hold exactly when $\Sigma$ is fixed set-wise by the set involution of $G\times G$ given by $(g,h) \mapsto (g^{-1},gh^{-1}g^{-1})$.

For the other condition, 

\begin{eqnarray*}
z\otimes z & = &\sum_{\substack{(g,h) \in \Sigma \\ (\gamma, \eta) \in \Sigma}} g\delta_h \otimes \gamma\delta_\eta  \\
\Delta(z) & = & \sum_{\substack{(g,h) \in \Sigma \\ s,t \in G \\ st = h}} g\delta_s \otimes g\delta_t \\
1\otimes z & = & \sum_{\substack{u \in G \\ (\gamma, \eta) \in \Sigma}} e\delta_u \otimes \gamma\delta_\eta
\end{eqnarray*}

Thus

\begin{eqnarray*}
(1\otimes z) \Delta(z) & = & \sum_{\substack{(g,h),(\gamma,\eta) \in \Sigma \\ s,t,u \in G \\ st=h}} g\delta_{g^{-1}sg}\delta_u \otimes \gamma g \delta_{g^{-1}\eta g}\delta_t \\
& = & \sum_{\substack{(g,h),(\gamma,\eta) \in \Sigma \\ s \in G \\ sg^{-1}\eta g=h}} g\delta_{g^{-1}sg} \otimes \gamma g \delta_{g^{-1}\eta g} \\ 
& = & \sum_{(g,h),(\gamma,\eta) \in \Sigma} g\delta_{g^{-1}hg^{-1}\eta^{-1}g^2} \otimes \gamma g \delta_{g^{-1}\eta g}.
\end{eqnarray*}

The requirement that this last sum be equal to $z\otimes z$ imposes more closure conditions on $\Sigma$:  for the second tensorand to range freely over all elements of $\Sigma$ without regard to the value of the first tensorand, $\Sigma$ must be closed under the action of the subgroup of $G$ generated by $\pi_1(\Sigma)$ acting on $G\times G$ by right multiplication on the first factor and right conjugation on the second.

Once this condition is met, letting $(\nu,\mu) = (\gamma g, g^{-1}\eta g)$ we can rewrite the expresion as

\[ \sum_{(g,h),(\nu,\mu) \in \Sigma} g\delta_{g^{-1}h\mu^{-1}g} \otimes \nu \delta_{\mu} \]

Thus for the first tensorand to range freely over $\Sigma$, $\pi_1^{-1}(g) = \{ k | (g,k) \in \Sigma\}$ must be fixed set-wise by the operations
$(g,h) \mapsto (g, g^{-1}h\mu^{-1}g)$ for all $\mu \in \pi_2(\Sigma)$.

These rather strange closure properties are easily seen to be satisfied of $\Sigma$ is chosen to be $H\times K$, for $K$ a finite subgroup of $G$ and $H$ a finite subgroup of its normalizer, thus giving us examples of traces with the necessary properties on a ribbon quasitriangular multiplier Hopf algebra which is not an ordinary ribbon Hopf algebra.

\section*{Acknowledgements}

 The work was partially supported by the National Natural Science Foundation of China (No. 11601231) 
 and the Fundamental Research Funds for the Central Universities (No. KJQN201716).  
 The first author wishes to thank the Kansas State University Department of Mathematics for hospitality 
 during his visit to the United States during which this work was carried out.


\vskip 0.6cm
\begin{thebibliography}{30}



 \bibitem{D10} Delvaux, L. (2010). Traces on multiplier Hopf algebras.
 \emph{Communications in Algebra}, 38(1): 346-360.

 \bibitem{DV04} Delvaux L. and Van Daele A. (2004). The Drinfel'd double versus the Heisenberg double for an algebraic quantum group.
 {\it Journal of Pure and Applied Algbera} 190: 59-84.

 \bibitem{DVW05} Delvaux, L., Van Daele, A. and Wang S. H. (2005). Quasitriangular (G-cograded) multiplier Hopf algebras.
 \emph{Journal of Algebra} 289: 484-514.

 \bibitem{DrVZ99} Drabant, B., Van Daele, A. and Zhang, Y. H. (1999). Actions of multiplier hopf algebras.
 \emph{Communications in Algebra} 27(9): 4117-4172.

 \bibitem{H96} Hennings, M. (1996). Invariants of links and 3-manifolds obtained from Hopf algebras.
 \emph{Journal of London Mathematical Society} 54(2): 594-624.
 
 \bibitem{I15} Ip, I. C. (2015). Positive representations of split real quantum groups: the universal $R$ operator.
 \emph{International Mathematics Research Notices} 1: 240-287.

 \bibitem{K95} Kassel, C. (1995). Quantum groups. Berlin: Springer.

 \bibitem{L10} Liu, G. H. (2010). Ribbon element on co-Frobenius quasitriangular Hopf algebras.
 \emph{Applied Mathematics} 1: 230-233.
 
 \bibitem{P11} Pennig U. (2011). Modular tensor categories from quantum doubles of finite groups.
  Lecture about "Topological Quantum Field Theories" during the winter semester 2011. 
  \url{http://www.math.uni-muenster.de/reine/u/upenn_01/slides/QuantumDouble.pdf}

 \bibitem{S} Shum, M.-C. (1994). Tortile tensor categories. 
 \emph{Journal of Pure and Applied Algebra} 93(1): 57-110.

 \bibitem{V94} Van Daele, A. (1994). Multiplier Hopf algebras.
 \emph{Transaction of the American Mathematical Society} 342(2): 917-932.

 \bibitem{V98} Van Daele, A. (1998). An algebraic framework for group duality.
 \emph{Advance in Mathematics} 140(2): 323-366.

 \bibitem{V08} Van Daele, A. (2008). Tools for working with multiplier Hopf algebras.
 \emph{The Arabian Journal for Science and Engineering} 33(2C): 505-527.

 \bibitem{YZM13} Yang, T., Zhou, X. and Ma T. S. (2013). On braided $T$-categories over multiplier Hopf algebras.
 \emph{Communications in Algebra} 41(8): 2852-2868.

 \bibitem{Z99} Zhang, Y. H. (1999). The quantum double of a coFrobenius Hopf algebra.
 \emph{Communications in Algebra}, 27(3): 1413-1427.




\end{thebibliography}
\end {document}